\documentclass[a4wide, oneside]{article}  

\usepackage{myPackages}
\usepackage{myAssignments3}
\usepackage{myOperators}

\usepackage{empheq}

\newcommand{\PD}{\operatorname{PD}}  

\newcommand{\D}{\Delta}

\usepackage{tikz}
\usetikzlibrary{arrows,decorations.pathmorphing,backgrounds,positioning,fit,petri}
\usetikzlibrary{calc,arrows.meta,positioning}
\usepackage{ragged2e}

\usepackage{fancyhdr,floatpag}
\floatpagestyle{fancy}

\usepackage{caption}
\usepackage{subcaption}

\DeclarePairedDelimiter\abs{\lvert}{\rvert}

\usepackage{fancyhdr,floatpag}
\floatpagestyle{fancy}

\usepackage{caption}
\usepackage{subcaption}

\usepackage{graphicx}

\titleformat{\section}
{\normalfont\normalsize\bfseries\center}{\thesection}{1em}{}

\titleformat{\subsection}
{\normalfont\normalsize\bfseries}{\thesection}{1em}{}

\begin{document}
\renewcommand{\abstractname}{\vspace{-\baselineskip}}

\title{Hypersurface singularities and virtually overtwisted contact structures}

\author{Edoardo Fossati}

\date{}

\newcommand{\Addresses}{{
  \bigskip

  E.~Fossati, \textsc{Scuola Normale Superiore, Piazza dei Cavalieri 7, Pisa.}\par\nopagebreak
  \textit{E-mail address:} \texttt{edoardo.fossati@sns.it}

}}

\maketitle

\begin{abstract}
It is known that the lens space $L(2n,1)$ supports a virtually overtwisted contact structure arising as the boundary of the Milnor fiber of a complex hypersurface singularity. In this article we study the problem of realizing other $(L(p,q),\xi)$ in such a way, obtaining a series of necessary conditions for this to happen. 
\end{abstract}

\section{Introduction}\label{intro}

Following \cite{nemethi}, we review some terminology of singularity theory. Let $f:(\C^3,0)\to (\C,0)$ be the germ of a complex analytic function and let 
\[K=f^{-1}(0)\cap S_{\varepsilon}\]
be the \emph{link of the singularity}, where $S_{\varepsilon}$ is the sphere of radius $\varepsilon$ centered at the origin.
Milnor proved in \cite{milnor} that there exists $\varepsilon_0>0$ such that $\forall\, 0<\varepsilon<\varepsilon_0$ the map
\[f/|f|:S_{\varepsilon}\smallsetminus K\to S^1=\{z\in\C:\,|z|=1\}\]
is a smooth fibration. For any such $\varepsilon$ there exists $\delta_{\varepsilon}$ so that $\forall\, 0<\delta<\delta_{\varepsilon}$ the restriction
\[
f:B_{\varepsilon}\cap f^{-1}(\partial D_{\delta}) \to \partial D_{\delta}
\]
is a smooth fibration, whose diffeomorphism type does not depend on $\varepsilon$ and $\delta$. This is what we refer to as the \emph{Milnor fibration} of $f$. The fiber 
\[F=F_{\varepsilon,\delta}=B_{\varepsilon}\cap f^{-1}(\delta)\]
is the \emph{Milnor fiber} of $f$ (we omit $\varepsilon$ and $\delta$ from the notation of $F$). If $f$ is the germ of an \emph{isolated} singularity, then we have a diffeomorphism $\partial F\simeq K$, but in the case of non-isolated singularity $K$ is not smooth (while $F$ and $\partial F$ are always smooth manifolds). The boundary of the Milnor fiber comes with an extra structure, as explained below. Recall that:

\begin{ndef}
A contact structure on a 3-manifold $M$ is a nowhere integrable planes distribution $\xi$. If there exists an embedded disk $D\hookrightarrow M$ such that $\xi$ agrees with $TD$ along the boundary $\partial D$, then the contact structure $\xi$ is said to be \emph{overtwisted}, otherwise it is called \emph{tight}. A tight structure $\xi$ is called \emph{universally tight} if its pullback to the universal cover $\widetilde{M}$ is tight. The tight structure $\xi$ is called \emph{virtually overtwisted} if its pullback to some finite cover is overtwisted.
\end{ndef}

\begin{nrem}
A consequence of the geometrization conjecture is that the fundamental group of any 3-manifold is \emph{residually finite} (i.e. any non trivial element is in the complement of a normal subgroup of finite index), and this implies that any tight contact structure is either universally tight or virtually overtwisted (see \cite{honda}).
\end{nrem}

Contact topology shows up in singularity theory in the following way: the Milnor fiber $F$ of a singularity comes with a Stein structure $J$ which makes it a Stein filling of its boundary $\partial F$ equipped with the contact structure 
\[\xi=T\partial F\cap JT\partial F,\]
which is always tight (see \cite{eliashberg}). Therefore, from a complex germ $f:(\C^3,0)\to (\C,0)$ we obtain a contact 3-manifold $(\partial F,\xi)$ with a Stein filling $(F,J)$ of it. We have a dichotomy: 

\begin{itemize}
\item the singularity is isolated. 
In this case the structure $\xi$ is universally tight (see \cite{lekili}). Another work on this topic is \cite{dynamics}, where the authors show that the link of the hypersurface singularity 
\[z^p+2xy=0\]
is $L(p,p-1)$ with its unique tight contact structure (universally tight). We will prove (see corollary \ref{L(p,p-1)}) that this is the only lens space arising as the link of an isolated hypersurface singularity.

\item The singularity is not isolated. By contrast with previous point, this is the only case where a virtually overtwisted contact structure can arise. A good source of examples is given by the \emph{Hirzebruch singularity}
\[z^2+xy^n=0\]
with $n>1$, for which the boundary of the associated Milnor fiber is $L(2n,1)$. This type of singularity is studied in \cite[section 6]{pichon}.
\end{itemize}

\begin{figure}[ht!]
\centering
\includegraphics[scale=0.35]{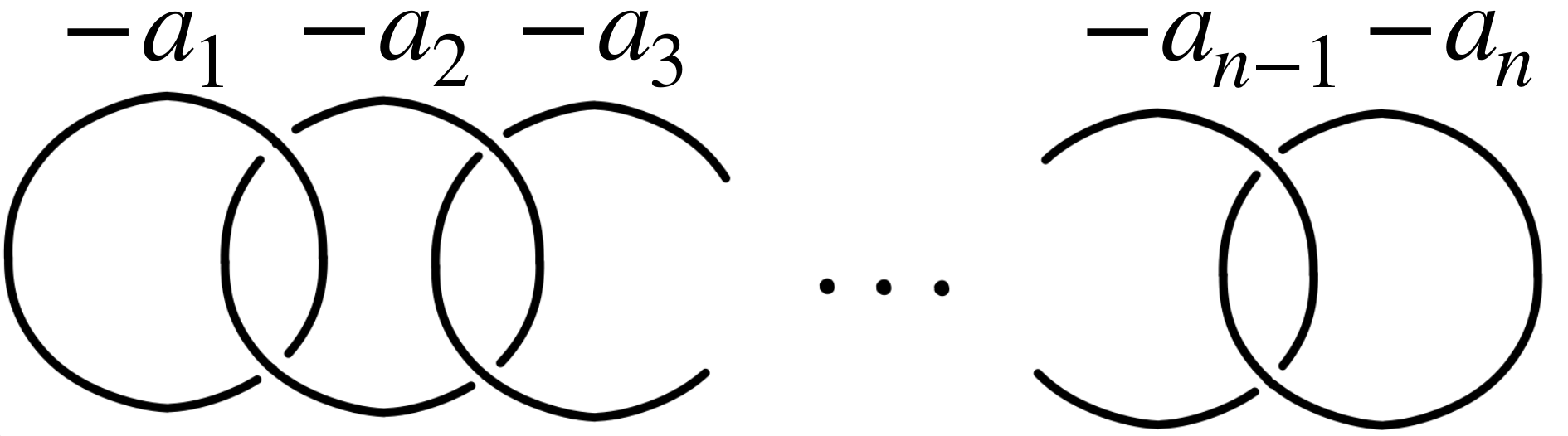}
\caption{Surgery link producing $L(p,q)$.}
\label{chain}
\end{figure}

\noindent Given a pair of coprime integers $p>q>1$, we consider the continued fraction expansion
\[\frac{p}{q}=[a_1,a_2,\ldots, a_n]=a_1-\frac{1}{a_2-\frac{1}{\ddots -\frac{1}{a_n}}}\]
with $a_i\geq 2$ for every $i$. As a smooth oriented 3-manifold, the lens space $L(p,q)$ is the integral surgery on a chain of unknots with framings $-a_1,-a_2,\ldots, -a_n$ (figure \ref{chain}). 

The goal of this article is to study those lens spaces $L(p,q)$ with a tight contact structure $\xi$ arising as the boundary of the Milnor fiber of a hypersurface singularity $f:(\C^3,0)\to (\C,0)$. Theorem \ref{mainthm} gives a partial answer to a question raised by \cite[open problems 24.4.2]{nemethi}.

\begin{thm}\label{mainthm}

Let $\xi_{vo}$ be a virtually overtwisted structure on $L(p,q)$. If we are in one of the cases below, then $(L(p,q),\xi_{vo})$ is not the boundary of the Milnor fiber of any complex hypersurface singularity:

\begin{itemize}

\item[a)] $p/q=[a_1,a_2,\ldots,a_n]$ and $a_i$ is odd for some $i$;

\item[b)] $p/q=[2x_1,2x_2]$;

\item[c)] $p/q=[2x_1,2x_2,\ldots,2x_n]$, with $x_i>1$ for every $i$ ($n\geq 3$) and either:
\begin{itemize}
\item[i)] $q^2\not\equiv 1 \mod p$;
\item[ii)] $q^2\equiv 1 \mod p$ and $n$ is even.
\end{itemize}

\end{itemize}
\end{thm}

The first step in proving this theorem is to characterize those contact structures $\xi$ which can appear in the context of hypersurface singularity: the fact that $c_1(\xi)$ vanishes imposes certain conditions on the coefficients in the continued fraction expansion of $p/q$, which allow us to prove part $(a)$ of theorem \ref{mainthm}. Section \ref{secvan} deals with these numerical restrictions, adapted to the language of contact geometry. To prove parts $(b)$ and $(c)$ we need to look closely at the topology of the Milnor fibration and analyze its monodromy. In order to derive our statements we study the integral orthogonal group of the intersection form of the Milnor fiber, imposing a further restriction coming from a theorem of \cite{acampo}. This is explained in section \ref{mainproof}. The article ends with an open question which focuses on the limits of theorem \ref{mainthm}: where are these techniques failing?

\paragraph{Acknowledgments} The author wishes to thank András Némethi for raising this problem and for the help in trying to tackle it. Thanks are due also to András Stipsicz for the hospitality at the Renyi Institute of Mathematics in Budapest, where this work originated. Lastly, a lot of support and key insights came from Paolo Lisca, who pointed out the results from \cite{gerstein} and \cite{c1planar}. This article will be part of the PhD thesis of the author.

\section{Vanishing of the rotation numbers}\label{secvan}

The goal of this section is to prove theorem \ref{thmc1=0}, which is a special case of \cite[corollary 1.5]{c1planar} (this result applies since, by \cite{schonenberger}, every contact structure on $L(p,q)$ is planar). We prove it using elementary techniques that do not involve the Heegaard-Floer contact invariant. Theorem \ref{thmc1=0} will be the starting point in the proof of theorem \ref{mainthm}. 

\begin{thm}\label{thmc1=0}
Let $L$ be a Legendrian linear chain of unknots in the standard contact $S^3$ and let $(L(p,q),\xi)$ be the contact 3-manifold obtained by Legendrian surgery on $L$. If $c_1(\xi)=0$, then $\rot (L_i)=0$ for every component $L_i$ of the link $L$.
\end{thm}

\begin{figure}[h!]
\centering
\begin{subfigure}[t]{.33\textwidth}
  \centering
  \includegraphics[scale=0.3]{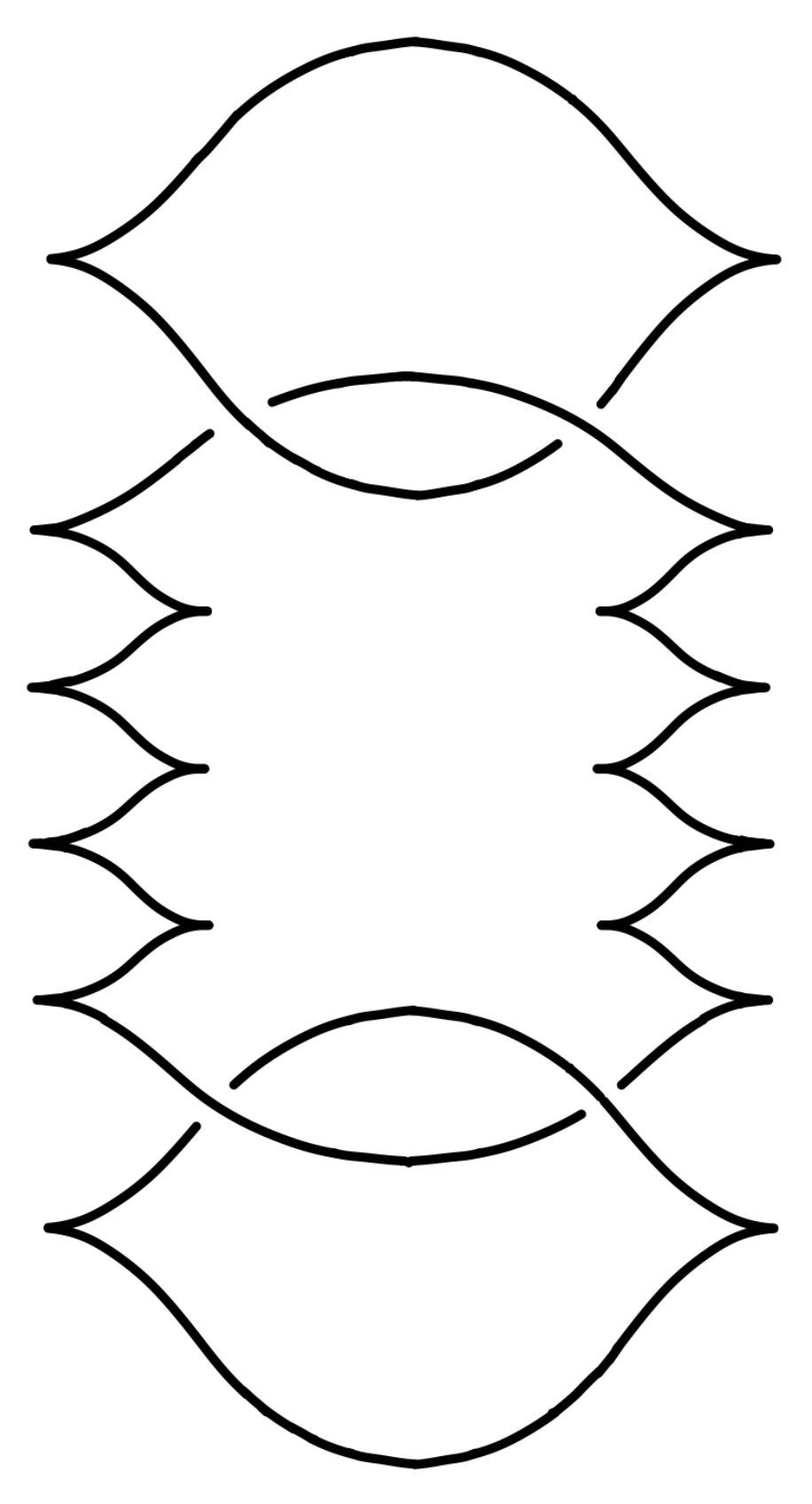}
  \caption{$L(28,15)$}
\end{subfigure}%
\begin{subfigure}[t]{.33\textwidth}
  \centering
 \includegraphics[scale=0.3]{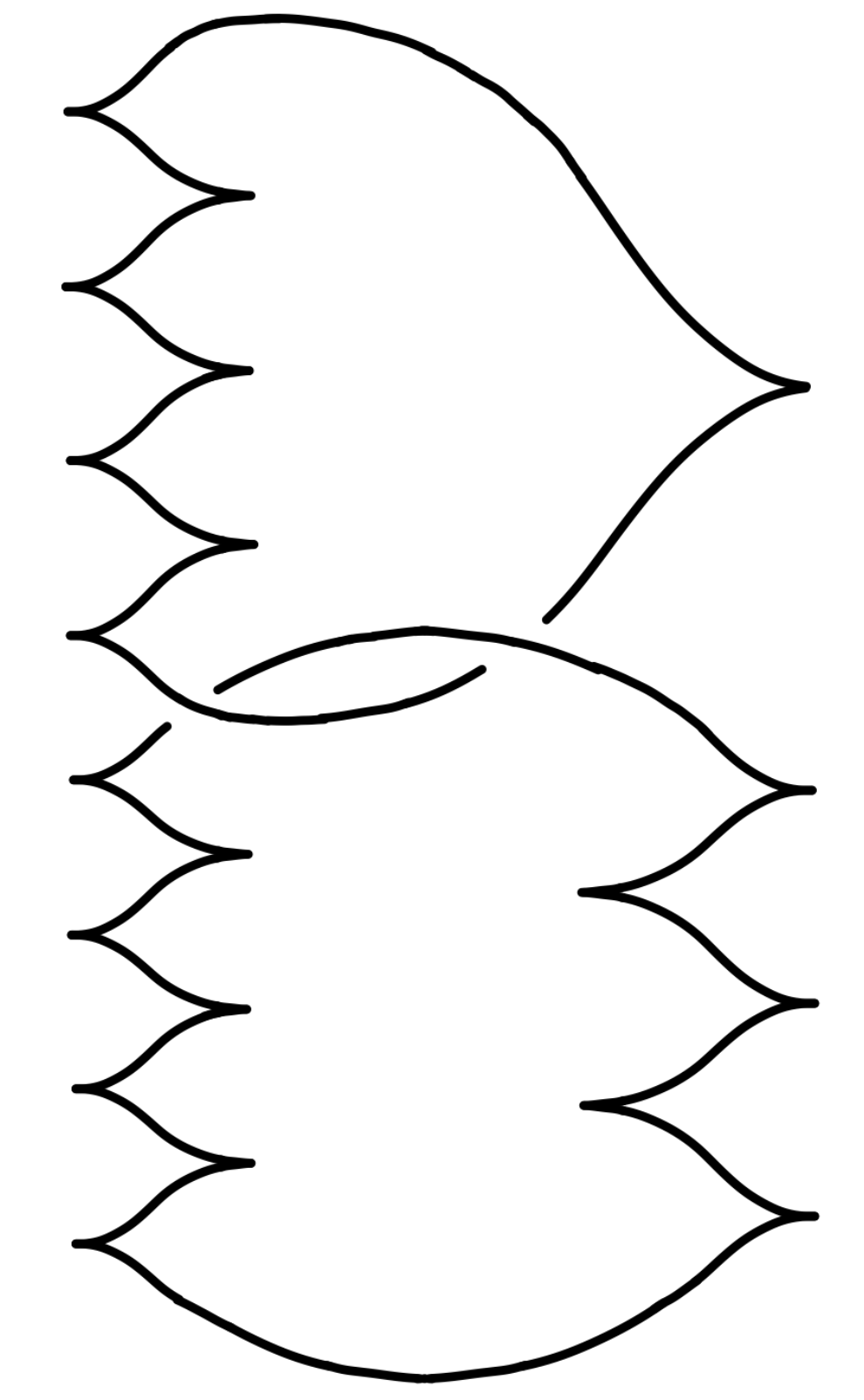}
 \caption{$L(34,7)$}
\end{subfigure}
\begin{subfigure}[t]{.33\textwidth}
  \centering
  \includegraphics[scale=0.3]{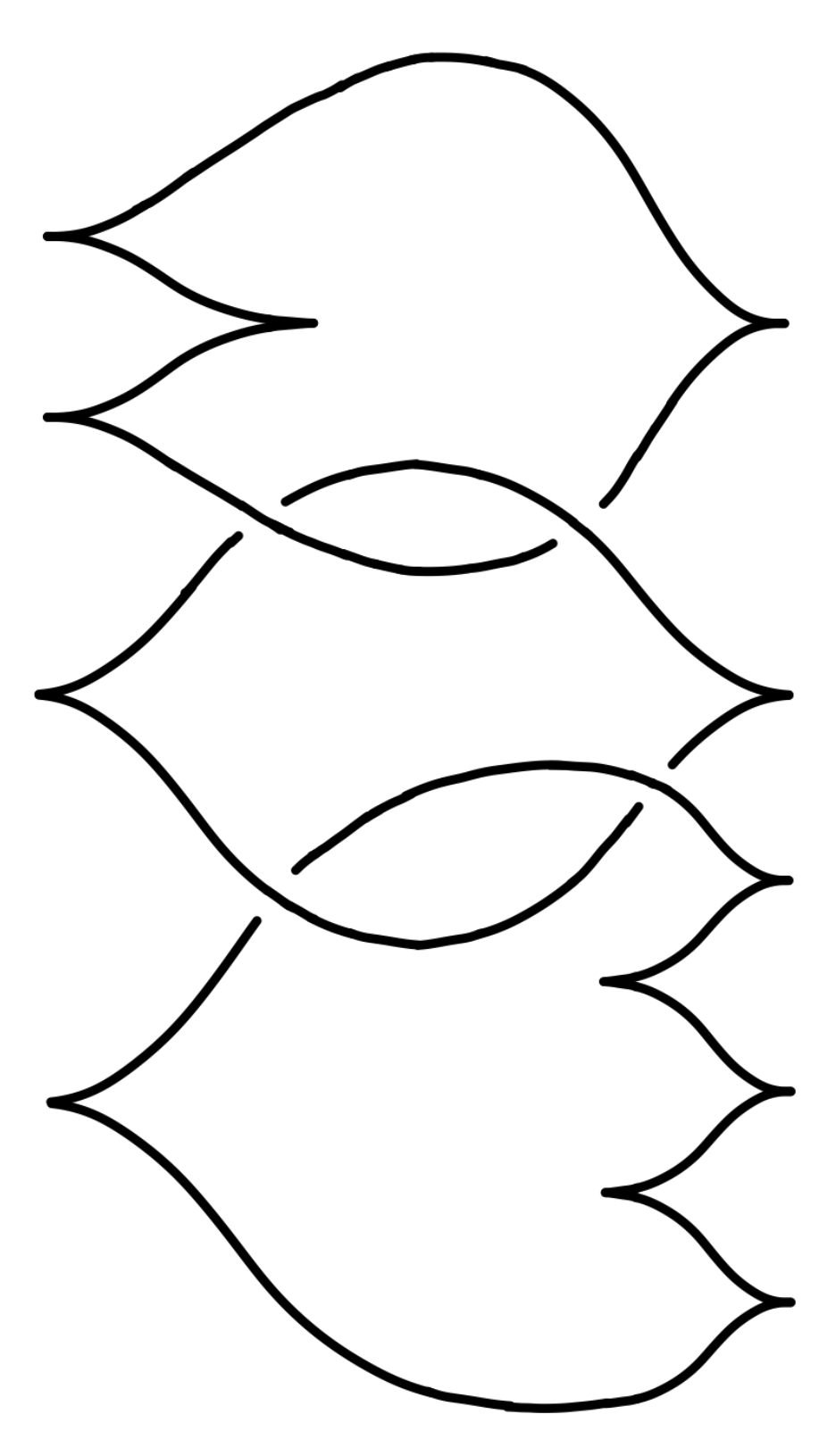}
  \caption{$L(17,7)$}
\end{subfigure}%
\caption{Examples of Legendrian links.}
\label{eg}
\end{figure}

To specify a tight contact structure on $L(p,q)$ we put the link of figure \ref{chain} into Legendrian position with respect to the standard tight contact structure of $S^3$, in order to form a linear chain of Legendrian unknots, see for example figure \ref{eg}. We do this in such a way that the Thurston-Bennequin number of the $i^{\mbox{th}}$ component is $a_i-1$. When such a Legendrian representative has been chosen, we get a tight structure $\xi$ on $L(p,q)$ by performing $(-1)$-contact surgery on it.

%

\noindent The rotation number of an \emph{oriented} Legendrian knot in $(S^3,\xi_{st})$ can be computed in the front projection by the formula
\[\rot(K)=\frac{1}{2}(c_D-c_U),\]
where $c_D$ and $c_U$ are respectively the number of down and up cusps. Since, as theorem \ref{thmc1=0} suggests, we are interested only in those Legendrian realizations with rotation number zero, we can avoid specifying an orientation.

To prove theorem \ref{thmc1=0} we need some notation that makes the computation easier, and a few results more. 
Let $Q$ be the linking matrix of the link of figure \ref{chain}, which represents the intersection form of the plumbed 4-manifold associated to the linear graph
\begin{center}
\begin{tikzpicture}
        \node[shape=circle,fill=black,inner sep=1.5pt,label=$-a_1$] (1)                  {};
        \node[shape=circle,fill=black,inner sep=1.5pt,label=$-a_2$] (2) [right=of 1] {}
        edge [-]               (1);
        \node[shape=circle,fill=black,inner sep=1.5pt,label=$-a_3$] (3) [right=of 2] {}
        edge [-]               (2);
        \node[shape=circle,fill=black,inner sep=1.5pt] (4) [right=of 3] {} ;
                 \node at ($(3)!.5!(4)$) {\ldots};

        \node[shape=circle,fill=black,inner sep=1.5pt,label=$-a_n$] (5) [right=of 4] {}
        edge [-]               (4);
\end{tikzpicture}
\end{center}
written in the ordered basis given by the meridian of each curve. We know that
\[H^2(L(p,q))\ni c_1(\xi)=0 \Leftrightarrow \PD(c_1(\xi))=0\in H_1(L(p,q)).\]
By looking at the linear plumbing graph we can express
\[\PD(c_1(\xi))=\sum_{i=1}^nr_i\mu_i,\]
where $r_i$ and $\mu_i$ are respectively the rotation number $\rot(L_i)$ and the meridian of the $i^{th}$ component (compare with \cite{ozbagci}). The $\mu_i$'s are related by a set of equations coming from $Q$:
\[
\begin{cases}
-a_1\mu_1+\mu_2=0	\\
\mu_1-a_2\mu_2+\mu_3=0 \\
\vdots \\
\mu_j-a_{j+1}\mu_{j+1}+\mu_{j+2}=0 \\
\vdots \\
\mu_{n-1}-a_n\mu_n=0.
\end{cases}
\]
By choosing $\mu_1\in H_1(L(p,q))$ as a generator we can consider it to be $1\in \Z/p\Z\simeq H_1(L(p,q))$ and therefore identify all the other $\mu_i$'s with numbers (which are well-defined modulo $p$). We derive relations which hold over $\Z$, but this will be enough to us thanks to lemma \ref{lemma0}.
Therefore we get recursive expressions of the form
\[
\begin{cases}
\mu_1=1\\
\mu_2=a_1\\
\mu_i=a_{i-1}\mu_{i-1}-\mu_{i-2}.
\end{cases}
\]
Then $c_1(\xi)$ is 0 exactly when 
\begin{equation}\label{chernrotations}
\sum_{i=1}^n r_i\mu_i\equiv 0\mbox{ modulo } \det(Q).
\end{equation}
Define
\[
\begin{cases}
\D[-1]=0;\\
\D[0]=1;\\
\D[i]=-a_i\D[i-1]-\D[i-2]
\end{cases}
\]
and note that
\begin{itemize}
\item $\det(Q)=\D[n]$;
\item $\sign(\D[i])=(-1)^i$, hence $\D[i]=(-1)^i|\D[i]|$;
\item $\mu_i> \mu_{i-1}$;
\item $|r_i|\leq a_i-2$ and $r_i\equiv a_i \mod 2$;
\item $\D[i]=(-1)^i\mu_{i+1}$ (proved by induction), therefore $|\D[i]|>|\D[i-1]|$.
\end{itemize}

\begin{lem}\label{lemma0}
Equation \ref{chernrotations} is satisfied in $\Z/\det(Q)\Z$ if and only if it is satisfied in $\Z$.
\end{lem}

\begin{prf}
We prove by induction on $n$ that 
\begin{equation}\label{inequality}
|\det(Q)|>\abs[\Big]{\sum_{i=1}^n r_i\mu_i}.
\end{equation}
This will tell that equation \ref{chernrotations} can only be satisfied with the zero-multiple of $\det(Q)$, since clearly for any non-zero $m\in\Z$ we have
\[|m\cdot \det(Q)|>|\det(Q)|.\]
If $n=1$ then the two sides of the inequality \ref{inequality} are respectively $|-a_1|$ and $|r_1|$, so it is true. 
The first interesting case is then $n=2$:
\[a_2a_1-1>r_1+r_2a_1?\] 
We have:
\begin{align*}
r_1+r_2a_2 < &(a_1-2)+(a_2-2)a_1\\
= & a_1a_2-a_1-2\\
< & a_1a_2-1.\qquad \checkmark
\end{align*}
The general case now:
\begin{align*}
\abs[\Big]{\sum_{i=1}^nr_i\mu_i}= & \abs[\Big]{\sum_{i=1}^nr_i\D[i-1]}\\
<&\left(\abs[\Big]{\sum_{i=1}^{n-1}r_i\D[i-1]}\right)+|r_n\D[n-1]| && \mbox{(induction)} \\
<&|\D[n-1]|+|r_n\D[n-1]|\\
\leq &|\D[n-1]|+|(a_n-2)\D[n-1]|.
\end{align*}
There are two possibilities for the right-hand side, according to the parity of $n$.
\begin{itemize}
\item[1)] $n$ is even (hence $\D[n]>0$):
\begin{align*}
|\D[n-1]+|(a_n-2)\D[n-1] |=& -\D[n-1]-a_n\D[n-1]+2\D[n-1]\\
=&-a_n\D[n-1]+\D[n-1]\\
<& -a_n\D[n-1]-\D[n-2] \mbox{ (because $|\D[i]|>|\D[i-1]|$ and $n$ is even)}\\
=&\D[n]\\
=&|\det(Q)|.
\end{align*}
\item[2)] $n$ is odd (hence $\D[n]<0$):
\begin{align*}
|\D[n-1]+|(a_n-2)\D[n-1] |=& a_n\D[n-1]-\D[n-1]\\
<& a_n\D[n-1]+\D[n-2] \mbox{ (because $|\D[i]|>|\D[i-1]|$ and $n$ is odd)}\\
=&-\D[n]\\
=&|\det(Q)|.
\end{align*}
\end{itemize}
\end{prf}

\begin{lem}\label{lemma1} 
\[a_n\mu_n-a_{n-1}\mu_{n-1}>0.\]
\end{lem}

\begin{prf} We prove it by induction on $n$. If $n=2$ then the formula is just
\[a_2\mu_2-a_1\mu_1=a_2a_1-a_1>0.\qquad\checkmark\]
In general, assuming the result true for $n-1$, we have
\begin{align*}
a_n\mu_n-a_{n-1}\mu_{n-1}=& a_n(a_{n-1}\mu_{n-1}-\mu_{n-2})-a_{n-1}\mu_{n-1}\\
=& a_{n-1}\mu_{n-1}(a_n-1)-a_n\mu_{n-2}&& \mbox{(induction)}\\
>& a_{n-2}\mu_{n-2}(a_n-1)-a_n\mu_{n-2}\\
=& (a_{n-2}a_n-a_{n-2}-a_n)\mu_{n-2}\\
>& 0.
\end{align*}
\end{prf}

\begin{lem}\label{lemma2} If $r_n\neq 0$, then we have
\[|r_n\mu_n|-|r_{n-1}\mu_{n-1}|>0.\]
\end{lem}

\begin{prf}
\begin{align*}
|r_n\mu_n|-|r_{n-1}\mu_{n-1}|>& |\mu_n|-(a_{n-1}-2)\mu_{n-1}\\
=& a_{n-1}\mu_{n-1}-\mu_{n-2}-a_{n-1}\mu_{n-1}+2\mu_{n-1}\\
=& 2\mu_{n-1}-\mu_{n-2}\\
>& 0.
\end{align*}
\end{prf}

\begin{lem}\label{lemma3}
If $r_n\neq 0$, then we have
\[|r_n\mu_n|-\abs[\Big]{\sum_{i=1}^{n-1}r_i\mu_i}>0.\]
\end{lem}

\begin{prf} We prove it by induction on $n$. If $n=2$ the formula is just
\[|r_2\mu_2|-|r_1\mu_1|=|r_2a_1|-|r_1|>a_1-(a_1-2)=2>0.\qquad\checkmark\]
Now we do the general case. Assume the inequality holds for $n-1$ and let $j\leq n-2$ be the biggest integer such that $r_j\neq 0$ (note that if $r_i=0,\;\forall i\leq n-2$, then by lemma \ref{lemma2} we would be done after the second line in the following computation). 
We have 
\begin{align*}
|r_n\mu_n|-\abs[\Big]{\sum_{i=1}^{n-1}r_i\mu_i} >& |\mu_n|-\abs[\Big]{\sum_{i=1}^{j-1}r_i\mu_i}-|r_{n-1}\mu_{n-1}| &&  \mbox{(induction)}\\
>& \mu_n-|r_j\mu_j|-|r_{n-1}\mu_{n-1}| \\
\geq & \mu_n-(a_j-2)\mu_j-(a_{n-1}-2)\mu_{n-1}\\
= & a_{n-1}\mu_{n-1}-\mu_{n-2}-(a_j-2)\mu_j-(a_{n-1}-2)\mu_{n-1} \\
= & 2\mu_{n-1}-\mu_{n-2}-(a_j-2)\mu_j	\\
> & \mu_{n-1}-(a_j-2)\mu_j\\
= & \mu_{n-1}+\mu_j+\mu_j-a_j\mu_j && (n-2\geq j\Rightarrow n-1\geq j+1)\\
>&\mu_{j+1}+\mu_{j-1}+\mu_{j-1}-a_j\mu_j && (\mu_{j+1}+\mu_{j-1}=a_{j}\mu_{j})\\
=&\mu_{j-1}\\
>&0.
\end{align*}
\end{prf}

\begin{lem} \label{lemma4}
\[\sum_{i=1}^{n}r_i\mu_i=0\Longrightarrow r_i=0,\;\forall i.\]
\end{lem}

\begin{prf}
\[\sum_{i=1}^{n}r_i\mu_i=0\Longrightarrow r_n\mu_n=-\sum_{i=1}^{n-1}r_i\mu_i\Longrightarrow |r_n\mu_n|=\abs[\Big]{\sum_{i=1}^{n-1}r_i\mu_i}.
\]
But if $r_n\neq 0$, then we should have a strict inequality by lemma \ref{lemma3}, hence $r_n=0$. By applying this repeatedly we get to
\[r_n=r_{n-1}=\ldots=r_1=0.\]
\end{prf}

\noindent We can finally give the following:

\begin{prf}[of theorem \ref{thmc1=0}]
By combining lemmas \ref{lemma0} and \ref{lemma4}, we have that
\[c_1(\xi)=0\Longleftrightarrow \sum_{i=1}^n r_i\mu_i\equiv 0\mod p \Longleftrightarrow \sum_{i=1}^n r_i\mu_i=0\Longleftrightarrow r_i=0\,\forall i,\]
where $r_i=\rot (L_i)$.
\end{prf}

\section{Proof of main theorem }\label{mainproof}

\begin{center}
Are the virtually overtwisted structures on lens spaces realizable as the boundary of the Milnor fiber of some complex hypersurface singularity?
\end{center}

\noindent In the universally tight case we know (see \cite{nemethippp} and \cite{bhupal}) that all the fillings come from algebraic geometry, and the article \cite{park} shows that complex surface singularities produce all the fillings of small Seifert 3-manifold equipped with the  canonical contact structure.

\vspace{0.5cm}

\noindent \textbf{Key fact:} suppose there is a polynomial function $f:(\C^3,0)\to (\C,0)$ such that 
\[\partial(F,J)=(L(p,q),\xi),\]
where $(F,J)$ is the Milnor fiber of $f$ and $\xi$ is the contact structure on the boundary induced by complex tangencies, as explained in section \ref{intro}. The links of figure \ref{eg} represent more than contact 3-manifolds: the components $L_i$ of any of those links can be thought as the attaching circles of the 2-handles of the Stein 4-manifold $(F,J)$. The first Chern class $c_1(F,J)\in H^2(F;\Z)$ evaluates on each 2-handle as the correspondent rotation number, and $c_1(\xi)$ is the restriction of $c_1(F,J)$. But the tangent bundle of $F$ is stably trivial, hence $c_1(F,J)=0$, and also $c_1(\xi)=0$ on the boundary $L(p,q)$. 

\begin{prf}[of theorem \ref{mainthm}a]
If $(L(p,q),\xi_{vo})$ is the boundary of the Milnor fiber of a complex hypersurface singularity, then $c_1(\xi_{vo})=0$ and, by theorem \ref{thmc1=0}, all $\rot(L_i)$ are zero. \\
\noindent Let $p/q=[a_1,a_2,\ldots,a_n]$ and remember that 
\[a_i\equiv \rot(L_i)\mod 2.\] 
Since all the rotation numbers are zero, the conclusion follows.
\end{prf}

Theorem \ref{thmc1=0} implies also the following corollary, which will be used later.

\begin{cor}\label{uniqueness}
If $(L(p,q),\xi)$ is the boundary of the Milnor fiber of a hypersurface singularity, then it has a unique Stein filling, which is the Milnor fiber itself.
\end{cor}

\begin{prf}
The fact that $\rot(L_i)=0$ for every $i$ implies, by \cite[theorem 1.3]{menke}, that from the chain of Legendrian unknots producing $(L(p,q),\xi)$ we can forget about those components with $\tb(L_i)\neq -1$ and look for Stein fillings of the 3-manifold $Y$ which is left. Then, all the Stein fillings of $(L(p,q),\xi)$ will be uniquely obtained by attaching the 2-handles (corresponding to the forgotten components) to the Stein fillings of $Y$. From the link diagram of $(L(p,q),\xi)$ we see that $Y$ is a connected sum of $(L(n_j,n_j-1),\xi_{st})$, where each of the prime factor corresponds to a string of $-2$ in the expansion of $p/q$. 
By \cite[theorem 16.9]{steinandback}, $Y$ admits a unique Stein filling, because each factor $(L(n_j,n_j-1),\xi_{st})$ does. This concludes the proof.
\end{prf}

\vspace{0.5cm}

\noindent Another consequence of theorem \ref{thmc1=0} is:

\begin{cor}\label{L(p,p-1)}
Let $(L(p,q),\xi)$ be a lens space with a contact structure arising as the link of an isolated hypersurface singularity. Then $q=p-1$.
\end{cor}

\begin{prf}
From theorem \ref{thmc1=0} we have that all the rotation numbers are zero and \cite[theorem 2.1]{lekili} says that $\xi$ is universally tight. By the work of \cite{honda} we know that a universally tight structure on a lens space is the result of contact $(-1)$-surgery on a link where all the stabilizations appear on the same side, i.e. when the (absolute values of the) rotation numbers are maximal. Therefore every Legendrian knot must have Thurston-Bennequin number equal to $-1$:
\[-\frac{p}{q}=[-2,-2,\ldots,-2]\qquad\Rightarrow\qquad q=p-1.\]
\end{prf}

\subsection*{Proof of Theorem 1b}

Corollary \ref{uniqueness} says that if $(L(p,q),\xi_{vo})$ arises as $\partial (F,J)$, then the Stein filling $F$ is uniquely determined: topologically it is given by the plumbing of spheres according to the expansion of $p/q$. The monodromy $\varphi$ of the Milnor fibration induces, in cohomology, a homomorphism
\[\varphi^*:H^*(F;\Z)\to H^*(F;\Z)\]
such that the alternating sum of the traces is zero, by \cite[theorem 1]{acampo}:
\[\trace(\varphi^0)-\trace(\varphi^1)+\trace(\varphi^2)-\trace(\varphi^3)+\trace(\varphi^4)=0.\]
In our case, the Stein fillings of those lens spaces with $c_1(\xi)=0$ are simply connected, hence $\varphi^1=\varphi^3=0$. Moreover, $\varphi^0:\Z\to\Z$ is the identity and $\varphi^4=0$. Therefore previous equation simply reads as:
\begin{equation}\label{trace}
1+\trace(\varphi^2)=0.
\end{equation}
The intersection form of $F$ must be preserved by the homomorphism $\varphi^2$ and we are therefore led to study its isometry group. If in this group there is no element whose trace is $-1$, then formula \ref{trace} cannot be satisfied, and we conclude that the polynomial function $f:(\C^3,0)\to (\C,0)$ with $(\partial F,\xi)=(L(p,q),\xi_{vo})$ does not exist.

\begin{prf}[of theorem \ref{mainthm}b]
Our goal is to prove that there is no $f:(\C^3,0)\to (\C,0)$, with non-isolated singularity, whose Milnor fiber has boundary $L(p,q)$. Note that 
\[x_1x_2>1,\]
otherwise the induced contact structure is universally tight. 

Assume by contradiction that such $f$ exists. Then, by corollary \ref{uniqueness}, we know that the Milnor fiber $F$ has negative-definite intersection form
\[-M=
\begin{bmatrix}
   -2x_1 & 1 \\
   1 & -2x_2 \\
\end{bmatrix}.
\]
By equality \ref{trace}, we must have $\trace(\varphi^2)=-1$. The morphism $\varphi^2:H^2(F;\Z)\to H^2(F;\Z)$ is induced by a diffeomorphism which preserves the intersection form and therefore it is represented by an integral matrix $A$ with
\[\begin{cases}
|\det(A)|=1 \\
\trace(A)=-1
\end{cases}
\]
and such that 
\[A(-M)A^T=-M.\]
We show now that such matrix cannot exist. We change sign to work with a positive definite matrix:
\[AMA^T=M\Longrightarrow
 \begin{bmatrix}
   a_1 & a_2 \\
   a_3 & a_4 \\
\end{bmatrix}
\begin{bmatrix}
   2x_1 & -1 \\
   -1 & 2x_2 \\
\end{bmatrix}
\begin{bmatrix}
   a_1 & a_3 \\
   a_2 & a_4 \\
\end{bmatrix} = 
\begin{bmatrix}
   2x_1 & -1 \\
   -1 & 2x_2 \\
\end{bmatrix}. \]
We get equations:
\[\begin{cases}
2x_1a_1^2-2a_1a_2+2x_2a_2^2=2x_1 \\
2x_1a_3^2-2a_3a_4+2x_2a_4^2=2x_2
\end{cases}
\]
that can be rewritten as

\begin{subequations}
\begin{empheq}[left=\empheqlbrace]{align}
  & (2x_1-1)a_1^2+(a_1-a_2)^2+(2x_2-1)a_2^2=2x_1 \label{star1}
  \\
  & (2x_1-1)a_3^2+(a_3-a_4)^2+(2x_2-1)a_4^2=2x_2. \label{star2} 
\end{empheq}
\end{subequations}

\noindent From equations \ref{star1} and \ref{star2} it follows that $a_1^2\leq 1$ and $a_4^2\leq 1$. But since $\trace(A)=a_1+a_4=-1$, we have that either
\[ \begin{cases}
a_1=0 \\
a_4=-1
\end{cases}
\qquad \mbox{or}\qquad 
\begin{cases}
a_1=-1 \\
a_4=0.
\end{cases}\]
We do the first case $(a_1=0,\,a_4=-1)$, the other one is the same. From $AMA^T=M$ we also get 
\[
 \begin{bmatrix}
   0 & a_2
\end{bmatrix}
\begin{bmatrix}
   2x_1 & -1 \\
   -1 & 2x_2 \\
\end{bmatrix}
\begin{bmatrix}
   a_3 \\
   -1  \\
\end{bmatrix} = -1, \]
which gives the equation $a_2(a_3+2x_2)=1$. Hence $a_2=a_3+2x_2=1$ (or both $-1$, but the conclusion is the same).
From \ref{star1} we have $x_2a_2^2=x_1$, which gives $x_2=x_1$. Then
\[\pm 1=\det(A)=\det 
\begin{bmatrix}
   0 & \pm 1 \\
   a_3 & -1  \\
\end{bmatrix} 
\Longrightarrow a_3=\pm 1.\]
We are in the case where $a_3+2x_2=1$, so either $x_1=x_2=0$ or $x_1=x_2=1$, which are both contradicting the condition $x_1x_2>1$.
\end{prf}

\subsection*{Proof of Theorem 1c}

Remember that in order to have $c_1(\xi)=0$ we need all the rotation numbers to be zero and in particular all the coefficients in the expansion to be even. Let $-M$ be the negative definite intersection lattice associated to the linear plumbing of spheres
\begin{center}
\begin{tikzpicture}
        \node[shape=circle,fill=black,inner sep=1.5pt,label=$-2x_1$] (1)                  {};
        \node[shape=circle,fill=black,inner sep=1.5pt,label=$-2x_2$] (2) [right=of 1] {}
        edge [-]               (1);
       \node[shape=circle,fill=black,inner sep=1.5pt,label=$-2x_3$] (3) [right=of 2] {}
        edge [-]               (2);  
        \node[shape=circle,fill=black,inner sep=1.5pt] (4) [right=of 3] {} ;
                         \node at ($(3)!.5!(4)$) {\ldots};

        \node[shape=circle,fill=black,inner sep=1.5pt,label=$-2x_n$] (5) [right=of 4] {}
        edge [-]               (4);
\end{tikzpicture}
\end{center}

We are looking for a matrix $A$ representing the monodromy of a Milnor fibration on the second cohomology group, that respect the intersection form of the Milnor fiber (i.e. $A(-M)A^T=M$) and whose trace is $-1$.

What we need to understand is the integral orthogonal group $O_{\Z}(-M)$ of the negative definite lattice $(\Z^n,-M)$, which is isomorphic to $O_{\Z}(M)$. The latter is studied in the article \cite{gerstein}, where the following theorem is proved:

\begin{thm*}[\cite{gerstein}] 
Let $M$ be the integer matrix
\[
\begin{bmatrix}
    2x_1       & -1  \\
 -1       & \ddots & \ddots \\
 & \ddots & \ddots & \ddots \\
  & & \ddots &\ddots &-1\\
  & & & -1 & 2x_n\\
\end{bmatrix}
\]
with $n\geq 2$ and $2x_i\geq 3\;\forall i$. 
\begin{itemize}
\item[i)] If $x_i\neq x_{n+1-i}$ for some $i$, then $O_{\Z}(M)=\{\pm \id\}$.
\item[ii)] If $x_i=x_{n+1-i}$ for every $i=1,\ldots, n$, then $O_{\Z}(M)=\{\pm \id,\pm \rho\}$, where $\rho$ is the isometry that inverts the order of a basis.
\end{itemize}
\end{thm*}

Therefore we can derive:

\begin{prf}[of theorem \ref{mainthm}c]
The condition $q^2\not\equiv 1 \mod p$ can be rephrased in terms of the coefficients of the expansion by saying that $x_i\equiv x_{n+1-i}$ for some $i$ (see \cite[appendix]{orlik}). On the other hand, $x_i=x_{n+1-i}$ for all $i$ if and only if $q^2= 1 \mod p$.

In the first case, previous theorem tells that if there is an integer matrix $A$ with $AMA^T=M$, then $A=\pm	\id$. Since this does not have trace $-1$, there cannot be a hypersurface singularity whose Milnor fiber has boundary $(L(p,q),\xi_{vo})$, otherwise we would have a contradiction with A'Campo's formula \ref{trace}.

In the second case, again A'Campo's formula cannot be satisfied because a hypersurface singularity would produce a Milnor fibration with monodromy $A\in\{\pm \id,\pm \rho\}$:
\[\rho=
\begin{bmatrix}
     & & & 1  \\
  &   & 1  \\
  & \reflectbox{$\ddots$} \\
   1\\
\end{bmatrix}
\]
and, if $n$ is even, then $\trace(\rho)=0\neq -1$.
\end{prf}

\begin{que} There are cases which are not covered by theorem \ref{mainthm}: every time that in the continued fraction expansion of $-p/q$ appears a $-2$, the orthogonal group $O_{\Z}(M)$ is harder to understand. Nevertheless, in the easier case when $q^2\equiv 1\mod p$ we have a complete description of $O_{\Z}(M)$, but inside this group there is a matrix with trace -1 if the length of the expansion is odd. A simple case is for example 
\[-\frac{p}{q}=-\frac{12}{7}=[-2,-4,-2].\] 
Our techniques indeed do not exclude that the isometry 
\[-\rho=
\begin{bmatrix}
    0 & 0 &-1  \\
  0 &-1  & 0  \\
   -1 &0  & 0  \\
\end{bmatrix}
\]
is the morphism induced by the monodromy of the Milnor fibration of a certain non-isolated hypersurface singularity producing $L(12,7)$. How can we deal with cases like this and exclude the existence of such a monodromy? Further works will hopefully clarify this problem and either find the polynomial function or rule out this possibility as well.
\end{que}

\bibliographystyle{alpha} 

\bibliography{Boundary_Milnor_fiber.bib}
\thispagestyle{plain}

\Addresses

\end{document}